\documentclass[11pt]{article}
\usepackage[latin1]{inputenc}
\usepackage{geometry}                
\usepackage{graphicx}
\usepackage{amssymb,amsmath}
\usepackage{epstopdf}
\usepackage{subfigure}

\usepackage{color,hyperref}
\hypersetup{
    colorlinks=true,
    linkcolor=blue,
    filecolor=magenta,      
    urlcolor=cyan,
}
 
\newcommand{\bx}{\mbox{\bf x}}
\newcommand{\bX}{\mbox{\bf X}}
\newcommand{\bD}{\mbox{\bf D}}
\newcommand{\by}{\mbox{\bf y}}

\newcommand{\bv}{\mbox{\bf v}}

\newcommand{\bu}{\mbox{\bf u}}

\newcommand{\FF}{\mbox{\bf F}}
\newcommand{\TT}{\mbox{\bf T}}
\newcommand{\ff}{\mbox{\bf f}}
\newcommand{\fg}{\mbox{\bf g}}

\newcommand{\fd}{\mbox{$\phi_{\epsilon}$}}
\newcommand{\dd}{\mbox{$\epsilon$}}

\title{Regularized Stokeslet segments}
\author{R. Cortez}
\date{}                                           

\begin{document}

\maketitle

\begin{abstract}
We present a variation of the method of regularized Stokeslet (MRS) specialized for the case of forces and torques distributed over filaments in three dimensions. The new formulation is based on the exact solution of Stokes equation generated by a linear continuous distribution of regularized forces along a line segment. Therefore, a straight filament with linearly varying forces does not require discretization. A general filament is approximated by a piecewise linear curve in three dimensions where the length of each line segment is chosen only based on the variation of the force field and the desired accuracy of its piecewise linear approximation. The most significant advantage of this formulation is that the values of the regularization parameter $\epsilon$ and the length of the segments $h$ are decoupled as long as $\epsilon<h$ so that $\epsilon$ can be selected as a proxy for the radius of the  filament and $h$ is chosen to discretize the forces and torques. We analyze the performance on test problems and present biological applications of sperm motility based on existing models of swimming flagella in open space and near a plane wall. The results show, for example, that because the forces along the flagellum vary mildly, a flagellum can be approximated with as few as 11 segments of length $h$ while fixing the regularization parameter to $\epsilon=h/30$, overcoming the need for hundreds of discretization nodes required by the MRS when $\epsilon$ is small.  The filament behaves like a slender cylindrical tube of radius $\approx 0.97\epsilon$ so that  the value of $\epsilon$ influences the flagellum's  swimming speed. For fixed regularization, doubling the number of line segments does not affect the results significantly as long as the force field is resolved.  Examples that require rotlets and potential dipoles along the filament are also presented.
\end{abstract}

\section{Introduction}
The method of regularized Stokeslets (MRS) is a popular method for computing viscous flows generated by external forces. For example, as a solid object moves through the fluid, its surface exerts a force causing fluid motion. The method has been used in biological applications ranging from sperm motility, microorganism swimming, cell motion, microfluidic devices, biofilm studies, and more. Part of the appeal of the method is its simplicity since the same formulation can be used for situations in which forces are distributed over surfaces, curves, or scattered points, leading to a unified expression of the form
\begin{equation}\label{MRS:generic}
\frac{d\hat{\bx}}{dt} = \bu(\hat{\bx}) = \sum_{j=1}^N S^\epsilon(\hat{\bx},\bx_j) \ff_j,
\end{equation}
where the force $\ff_j$ is exerted at $\bx_j$ and $\hat{\bx}$ is the evaluation point. 
This equation can be used to evaluate the velocity field at some locations in the fluid or at the same points where the forces are applied.
The kernel $S^\epsilon$ is derived by finding the
exact solution of the Stokes equations, $0=-\nabla p + \mu \Delta\bu + \ff \fd$, $\nabla\cdot\bu=0$
in $\mathbb R^3$, when $\fd$ is a smooth function that approximates the
Dirac delta distribution. Here we consider the spherically symmetric cutoff function $\fd(\bx) = 15 \epsilon^4/8\pi(|\bx|^2+\epsilon^2)^{7/2}$, where $\epsilon$ is a small parameter to be chosen~\cite{Cortez:01, Cortez:2005}.

Without regularization, a force is of the form $\ff \delta(\bx-\bx_j)$ and the resulting kernel is the
Stokeslet
\[
8\pi\mu S^0(\hat{\bx},\bx_j) = \frac{I}{\|\hat{\bx}-\bx_j\|} + \frac{(\hat{\bx}-\bx_j)(\hat{\bx}-\bx_j)^T}{\|\hat{\bx}-\bx_j\|^3}
\]
which is singular at the locations where forces are applied. The regularization of the kernel through the use of smooth approximations of the delta distribution eliminates the singularity.
The difference in each of the cases mentioned above is the interpretation of the terms in the sum. 
For scattered points the summation represents a superposition of contributions from regularized point forces.  
For the case of forces on a surface, the summation represents a discrete version of a surface integral where each term is the integral over a small patch on the surface. In this case, the regularization could be removed completely resulting in a principal value integral that
can be computed in other ways~\cite{Pozrikidis:92}.  Nevertheless, the regularization can be helpful as a numerical treatment of a weakly singular integral, particularly in combination with  finite elements~\cite{SmithBEM, SMITHNN, Montenegro:2016}, discretization refinement~\cite{BARREROGIL} or error reduction formulas~\cite{nguyencortez2014}. These approaches tend to
reduce the sensitivity of the velocity field to the regularization parameter. 

There are important applications involving the simulation of cilia and flagella
that are modeled as curves in three-dimensional space with a force distribution on it.
The flagella may be associated with bacteria or spermatozoa. For instance, the study of sperm
motility includes reaching a full understanding of the many waveforms observed
experimentally, the characterization of swimming trajectories, and chemical considerations.
The method of regularized Stokeslets has been useful in sperm motility for the study of
hyperactivated waveforms~\cite{Olson:2011}, interaction with a planar 
wall~\cite{CCDGK, spagnolie2012, JulieEtAl},
sperm-sperm interactions~\cite{SIMONS20151639}, 
 bundling of bacterial flagella~\cite{Flores:2005, Cisneros2008},
 and more~\cite{Gilles:2009, Ishimoto2017, wrobel2016enhanced}.
From a mathematical standpoint, the summation in Eq.~\eqref{MRS:generic}
for forces along a curve represents a line integral of the Stokeslet, which is divergent when evaluated
on the curve without regularization. 
The MRS addresses the singularity by spreading each force over a small sphere through the function $\fd$.
In computations, the curves are discretized using nodes separated by a distance $h$ with forces at the nodes. Since the regularization parameter $\epsilon$ controls the width of the function $\fd$, in practice $\epsilon$ is chosen large enough for contiguous functions to overlap but small enough for $\fd$ to be a reasonable approximation of the delta function. If $\epsilon$ is too small compared to $h$, contiguous cutoff functions do not overlap sufficiently and fluid leaks through the filament.  Thus we may adjust the regularization parameter to keep $\epsilon$ as small as possible while requiring the leak to be within a tolerance.

In scientific applications sometimes the regularization parameter is chosen based on physical arguments. For instance, in sperm motility $\epsilon$ may be chosen based on the sperm flagellum radius. In this case the discretization size $h$ is adjusted based on the value of $\epsilon$. In practical computations the regularization parameter is proportional to $h$ with proportionality constant typically between 1/2 and 7, depending on the cutoff being used. 
Representative dimensions of mammalian sperm result in a ratio of flagellum length to radius of 
about 100~\cite{Ishimoto20150172}.  Sea urchin sperm are about 40-50 $\mu$m long and  0.1-0.2 $\mu$m
in diameter~\cite{Kinukawa02022007}, for a length-to-radius ratio of at least 400. This means that setting
$\epsilon$ to be the flagellum radius and $h\approx\epsilon$ would require a discretization of about 400 nodes.

This article proposes a new way of using regularized Stokeslets for applications in which the external forces are distributed along a curve in three dimensional space. The goal is to develop a framework in which one can use smaller values of $\epsilon$ than currently done without requiring a large number of discretization nodes. Figure 1 shows results from a simulation of a swimming sperm with $\epsilon= 0.005$ in which a flagellum of length one has been represented using only 11 straight segments.
The idea is to approximate the given curve as connected line segments, along which the force density is assumed to vary linearly. The fluid velocity then includes the contribution from a continuous force distribution along each segment. The method is derived in Section 2, including a  recursion formula that is useful for the implementation of the method. Section 3 describes numerical examples and the implementation of the no-flow boundary condition in the plane $z=0$ for simulations of
flows in the half space.

\section{Elements of Regularized Stokes Flow}
For the particular choice of regularization~\cite{Cortez:2005} $\phi(R_0) = 15 \epsilon^4  R_0^{-7}/8\pi$, with $R_0^2 =  |\bx|^2+\epsilon^2$, the
velocity at the evaluation point $\hat{\bx}$ due to a regularized force $\ff$ applied at a point $\by_0$ is the regularized Stokeslet
\begin{equation}\label{regStokeslet}
8\pi\mu \bu(\hat{\bx}) =  \left( \frac{1}{R_0}+\frac{\epsilon^2}{R_0^3} \right) \ff +  \frac{(\ff\cdot\bx)\bx}{R_0^3}
\end{equation}
where $\bx=\hat\bx-\by_0$. 
Other solutions are derived by differentiation. For example, the fluid velocity due to a torque $\bf\tau$ is
the curl of Eq.~\eqref{regStokeslet}
\begin{equation}\label{regRotlet}
8\pi\mu \bu =  \left( \frac{2}{R_0^3}+\frac{3\epsilon^2}{R_0^5} \right)({\bf\tau}\times\bx)
\end{equation}
and a potential dipole of strength $\fg$ is
\begin{equation}\label{regDipole1}
8\pi\mu \bu =  -\left( \frac{2}{R_0^3}-\frac{6\epsilon^2}{R_0^5} \right) \fg +  \frac{6 (\fg\cdot\bx)\bx}{R_0^5}.
\end{equation}

A different regularization of the dipole of strength $\fg$, used in the Kirchhoff rod model~\cite{OLSON2013169}, is
\begin{equation}\label{regDipole2}
8\pi\mu \bu =  -\left( \frac{2}{R_0^3}+\frac{3\epsilon^2}{R_0^5} - \frac{15\epsilon^4}{R_0^7} \right) \fg 
+ (\fg\cdot\bx)\bx \left( \frac{6 }{R_0^5} + \frac{15\epsilon^2}{R_0^7} \right).
\end{equation}

\subsection{Flow due to a segment of Stokeslets}

\begin{figure}[hbtp]
\begin{picture}(200,50)
\put(20,10){\circle*{5}}
  \put(20,10){\vector(1,2){7}}  
\put(90,10){\circle*{5}}
  \put(90,10){\vector(3,1){22}}  
\put(60,10){\circle*{5}}\put(50,0){$\by(\alpha)$}
\put(13,0){$\by_0$}\put(83,0){$\by_1$}
\put(20,27){$\ff_0$}\put(106,22){$\ff_1$}
  \put(53,45){\circle*{3}}   \put(50,48){$\hat{\bx}$} 
   \put(20,10){\line(1,0){70}}
   \put(160,40){Notation: $\bv = \by_0-\by_1$,\ \ $|\bv|=\ell$}
   \put(210,20){$\bx = \hat{\bx}-\by(\alpha)$,\ \ $\bx_1=\hat{\bx}-\by_1$,\ \ $\bx_0=\hat{\bx}-\by_0$}
   \put(210,0){$R^2 = |\bx|^2+\dd^2$, $R_1^2 = |\bx_1|^2+\dd^2$, $R_0^2 = |\bx_0|^2+\dd^2$}
  \end{picture}
  \caption{Notation and schematic of a line segment from $\by_0$ to $\by_1$. The segment
  is along the vector $\bv$ and the force density along the segment is linearly interpolated from the
  end points. The velocity is evaluated at the point $\hat\bx$.}
  \label{fig:schematic}
  \end{figure}

We consider a line segment $\by(\alpha)=\by_0 + \alpha(\by_1- \by_0)$ of length
$\ell$ between the endpoints $\by_0$ and $\by_1$.  The parameter $\alpha\in[0,1]$ is dimensionless
and we assume a linear force density
$\ff(\alpha) = \ff_0 + \alpha(\ff_1- \ff_0)$ (see Figure~\ref{fig:schematic}).  
The force on an infinitesimal section of the segment is $ \ff(\alpha) \ell d\alpha$ and the net force is 
\[
F   = \ell \int_0^1 \ff(\alpha)  d\alpha  = \frac{\ell}{2}(\ff_0+\ff_1).
\]

The goal is to derive a formula for the
fluid velocity due to the continuum of forces. 
The velocity at ${\hat\bx}$ due to $\ff=\ff_a+ \alpha \ff_b$ along the straight segment
$\by=\by_0-\alpha \bv$ is 
\[
8\pi\mu \bu = \ell \int_0^1 \left( \frac{1}{R}+\frac{\epsilon^2}{R^3} \right) \ff +\frac{(\ff\cdot\bx)\bx}{R^3} d\alpha
\]
where $\bx={\hat\bx}-\by(\alpha)$ and $R^2=|\bx|^2+\epsilon^2=|\hat\bx-\by(\alpha)|^2+\epsilon^2$.  
The flow due to other elements is defined similarly. 
The important feature to notice is that all of these formulas consist of terms of the form
$P(\alpha) R_0^{q}$
where $P(x)$ is a polynomial with vector coefficients and $q$ is a nonzero integer. 
We will take advantage of this structure in the next sections.

\subsection{A Continuous Load on a Line Segment}
We begin with a straight line segment
$\by(\alpha)=\by_0-\alpha \bv$ where $\bv$ is a constant vector and 
$\bx=\hat{\bx}-\by(\alpha)$ with $R^2=|\bx|^2+\epsilon^2$.
If $\ff(\alpha)$ is a polynomial in $\alpha$, then so is $(\ff\cdot\bx)\bx$ so that the velocity field 
generated by this element is a linear combination of terms of the form
\begin{equation}\label{Tnq_definition}
T_{n,q} = \int_0^1 \alpha^n R(\alpha)^q d\alpha.
\end{equation}

The following identities can be established for a line segment of length $L$ and $m\ge0$:
\begin{enumerate}
\item $\bx\cdot\bv = \bx_0\cdot\bv+\alpha L^2$
\item $\bx'(\alpha)=-\by'(\alpha)=\bv$
\item $
L^2 R(\alpha)^2 - (\bx(\alpha)\cdot\bv)^2 = 
L^2 R_0^2 - (\bx_0\cdot\bv)^2\ \ \forall \alpha\in[0,1] $
\item $ \frac{d}{d\alpha} \left(  \alpha^m  \frac{R^p}{L^2}  \right)  =
 \frac{m}{L^2}  \alpha^{m-1} R^p
+p\alpha^{m+1} R^{p-2} + \frac{p(\bx_0\cdot\bv)}{L^2}\alpha^m  R^{p-2}$,\ \ for $m\ge0$
\end{enumerate}

Integrating the last formula and using Eq.~\eqref{Tnq_definition} we arrive at the recursion for $q\ne-2$
\begin{equation}\label{Recurrence}
T_{n,q} =
  \frac{ \alpha^{n-1} R^{q+2}}{(q+2)L^2} \Bigg|_0^1 
- \frac{n-1}{(q+2)L^2} T_{n-2,q+2} 
- \frac{(\bx_0\cdot\bv)}{L^2}  T_{n-1,q}, \ \ \ \ \ \  n\ge1
\end{equation}
This is similar to the recurrence formula derived by Chwang and  Wu in their derivation of the
flow past a prolate spheroid~\cite{ChwangWu74}.
Below are specific examples for forcing coefficients that vary linearly along the segment.

\subsection{The regularized Stokeslet on a line segment}
The velocity at $\hat{\bx}$ due to force density $\ff=\ff_a+ \alpha \ff_b$ along the  line segment
$\by(\alpha)=\by_0-\alpha \bv$ is
\begin{equation}\label{def:Stokeslet}
8\pi\mu \bu(\hat{\bx}) = L \int_0^1 \left( \frac{1}{R}+\frac{\epsilon^2}{R^3} \right) \ff +  \frac{(\ff\cdot\bx)\bx}{R^3} d\alpha
\end{equation}
where $\bx=\hat{\bx}-\by(\alpha)$ and $R^2=|\bx|^2+\epsilon^2$.
 The velocity can be written as
\[
(8\pi\mu/L) \bu(\hat{\bx}) =  \ff_a (T_{0,-1}+\epsilon^2 T_{0,-3})+\ff_b(T_{1,-1}+\epsilon^2 T_{1,-3})+\sum_{n=0}^3 \ff_n T_{n,-3}
\]
where the coefficients $\ff_n$ are  
\begin{eqnarray*}
\ff_0 &=& (\ff_a\cdot\bx_0)\bx_0, \\
\ff_1 &=& (\ff_a\cdot\bv)\bx_0 +(\ff_a\cdot\bx_0)\bv + (\ff_b\cdot\bx_0)\bx_0 ,\\ 
\ff_2 &=& ( (\ff_a\cdot\bv)\bv + (\ff_b\cdot\bx_0) )\bv + (\ff_b\cdot\bv)\bx_0, \\ 
\ff_3 &=&  (\ff_b\cdot\bv)\bv.
\end{eqnarray*}

Using the recursion~\eqref{Recurrence} with $q=-3$ we have
\[
T_{n,-3} =
 - \frac{ \alpha^{n-1} R^{-1}}{L^2} \Bigg|_0^1 
+ \frac{n-1}{L^2} T_{n-2,-1} 
- \frac{(\bx_0\cdot\bv)}{L^2}  T_{n-1,-3}
\]
beginning with the values when $n=0$ computed by direct integration, we get
the sequence of terms required for the segment of Stokeslets
\begin{eqnarray*}
T_{0,1} &=&  \frac{R (\bx\cdot\bv)}{2 L^2} + \frac{ (L^2R_0^2-(\bx_0\cdot\bv)^2) }{2 L^3} \log(LR+(\bx\cdot\bv)) \Bigg|_0^1  \\
T_{0,-1} &=& \frac{1}{L}\log(LR+(\bx\cdot\bv))  \Bigg|_0^1   \\
T_{0,-3} &=&  \frac{ -1}{R(LR+(\bx\cdot\bv))} \Bigg|_0^1   \\[8pt]
T_{1,-1} &=&  \frac{ R}{L^2}\Bigg|_0^1   - \frac{(\bx_0\cdot\bv)}{L^2}  T_{0,-1} \\
T_{1,-3} &=&
 - \frac{  R^{-1}}{L^2}  \Bigg|_0^1 
- \frac{(\bx_0\cdot\bv)}{L^2}  T_{0,-3}   \\
T_{2,-3} &=&
 - \frac{ \alpha R^{-1}}{L^2}  \Bigg|_0^1 
+ \frac{1}{L^2} T_{0,-1} 
- \frac{(\bx_0\cdot\bv)}{L^2}  T_{1,-3}   \\
T_{3,-3} &=&
 - \frac{ \alpha^{2} R^{-1}}{L^2}  \Bigg|_0^1 
+ \frac{2}{L^2} T_{1,-1} 
- \frac{(\bx_0\cdot\bv)}{L^2}  T_{2,-3}
\end{eqnarray*}

In our implementation it seems sufficient to use these terms as shown here. Alternatively,
one can use nested substitution to write a final formula only in terms of
$T_{0,p}$ for a few values of $p$.

\subsection{The regularized Dipole}
Similarly, the velocity at $\hat{\bx}$ due to a dipole $\fg=\fg_a+ \alpha \fg_b$ along the straight segment
$\by=\by_0-\alpha \bv$ is
\begin{equation}\label{def:dipole}
8\pi\mu \bu = L \int_0^1 -\left( \frac{2}{R^3}-\frac{6\epsilon^2}{R^5} \right) \fg +  \frac{6 (\fg\cdot\bx)\bx}{R^5} d\alpha
\end{equation}
where $\bx=\hat{\bx}-\by(\alpha)$ and $R^2=|\bx|^2+\epsilon^2$.
The velocity can be written as
\[
(8\pi\mu/L) \bu = -\fg_a (2T_{0,-3}-6\epsilon^2 T_{0,-5}) - \fg_b(2T_{1,-3}-6\epsilon^2 T_{1,-5})+6\sum_{n=0}^3 \fg_n T_{n,-5}
\]
where the coefficients $\fg_n$ are the same as $\ff_n$ in the Stokeslet case. Using the recursion~\eqref{Recurrence} 
with $q=-5$ gives
\[
T_{n,-5} =
  \frac{- \alpha^{n-1} R^{-3}}{3L^2} \Bigg|_0^1 
+ \frac{n-1}{3L^2} T_{n-2,-3} 
- \frac{(\bx_0\cdot\bv)}{L^2}  T_{n-1,-5},\ \ \ \ \ \ n\ge1.
\]
The required terms are
$T_{0,-1}$, $T_{1,-1}$
$T_{0,-3}$, $T_{1,-3}$, 
$T_{0,-5}$, 
$T_{1,-5}$, 
$T_{2,-5}$, and 
$T_{3,-5}$.

\subsection{The regularized Rotlet}
The final example is the velocity at ${\hat\bx}$ due to a rotlet $\tau=\tau_a+ \alpha \tau_b$ along the straight segment
$\by=\by_0-\alpha \bv$ 
\[
8\pi\mu \bu = L \int_0^1 \left( \frac{2}{R^3}+\frac{3\epsilon^2}{R^5} \right) \ (\tau\times\bx)\ d\alpha
\]
where $\bx={\hat\bx}-\by(\alpha)$ and $R^2=|\bx|^2+\epsilon^2$.
The velocity can be written as
\[
(8\pi\mu/L) \bu = \sum_{n=0}^2 \tau_n (2T_{n,-3}+3\epsilon^2 T_{n,-5})
\]
where the coefficients $\tau_n$ are 
$
\tau_0 = (\tau_a\times\bx_0),\ \ \ 
\tau_1 = (\tau_b\times\bx_0)+(\tau_a\times\bv),\ \ \ 
\tau_2 = (\tau_b\times\bv).
$
Using the recursion, we will need
$T_{0,-1}$, $T_{1,-1}$
$T_{0,-3}$, $T_{1,-3}$, $T_{2,-3}$,  
$T_{0,-5}$, 
$T_{1,-5}$, and
$T_{2,-5}$.

\subsection{Enforcing a prescribed velocity on a segment}
The velocity at a point $\hat\bx$ due to force density along a curve 
of length $\ell$ is
\begin{equation}\label{def:ContOperator}
8\pi\mu \bu =  \int_0^\ell \left( \frac{1}{R}+\frac{\epsilon^2}{R^3} \right) \ff_s +  \frac{(\ff_s\cdot\bx)\bx}{R^3} ds
\end{equation}
where $s$ is the arclength
can be inverted to find a force density that satisfy a specified velocity boundary condition.
We discretize the curve using $N_n$ nodes ($N_n-1$ segments) denoted by $\{ {\by}_k\}_{k=1}^{N_n}$ and consider a
linear force density on each segment, $\ff_s=\ff_{k+1}+\alpha(\ff_{k+1}-\ff_k)$. Then the velocity at a point $\hat\bx$ can be written as
\begin{equation}\label{def:DiscOperator}
8\pi\mu \bu(\hat\bx) = \sum_{k=1}^{N_n -1}  L_k \left(  \ff_k (T^{k,k+1}_{0,-1}+\epsilon^2 T^{k,k+1}_{0,-3})
                               +(\ff_{k+1}-\ff_k)(T^{k,k+1}_{1,-1}+\epsilon^2 T^{k,k+1}_{1,-3})        
                               +\sum_{n=0}^3 \ff_{n,k} T^{k,k+1}_{n,-3} \right)
\end{equation}
where we use superscript $k,k+1$ to emphasize that the terms are for the segment between $\by_k$ and $\by_{k+1}$.
Denoting $\by_k-\by_{k+1}$ by $\bv_k$ and $L_k = \|\bv_k\|$, the coefficients $\ff_{n,k}$ are  
\begin{eqnarray*}
\ff_{0,k} &=&  (\by_k\by_k^T) \ff_k, \\
\ff_{1,k}  &=& (\by_k\bv_k^T+\bv_k\by_k^T)\ff_k + (\by_k\by_k^T)(\ff_{k+1}-\ff_k) \\
\ff_{2,k}  &=& (\bv_k\bv_k^T)\ff_k +(\bv_k\by_k^T + \by_k\bv_k^T)(\ff_{k+1}-\ff_k) \\
\ff_{3,k}  &=&  (\bv_k\bv_k^T) (\ff_{k+1}-\ff_k).
\end{eqnarray*}
 Then $8\pi\mu \bu(\hat\bx)  =  M_1 \ff_k + M_2 \ff_{k+1}$ where the $3\times3$ blocks are
{\small
\begin{eqnarray}
M_2 &=& L_k \left(    (T^{k,k+1}_{1,-1}+\epsilon^2 T^{k,k+1}_{1,-3}) 
+ T^{k,k+1}_{1,-3}  (\by_k\by_k^T)
+ T^{k,k+1}_{2,-3}  (\by_k\bv_k^T+\bv_k\by_k^T)
+ T^{k,k+1}_{3,-3} (\bv_k\bv_k^T)  \right)        \label{eq:M1}\\
M_1 &=& L_k \left(  (T^{k,k+1}_{0,-1}+\epsilon^2 T^{k,k+1}_{0,-3}) I 
+ T^{k,k+1}_{0,-3}(\by_k\by_k^T)
+ T^{k,k+1}_{1,-3} (\by_k\bv_k^T+\bv_k\by_k^T)
+ T^{k,k+1}_{2,-3}(\bv_k\bv_k^T) \right)
- M_2.      \label{eq:M2}
\end{eqnarray}
}

We use these expressions to build a larger matrix by evaluating the expression in~\eqref{def:DiscOperator} at
each of the nodes. The final linear system describes the 
relation between the velocity of nodes describing a piecewise linear curve
to the force density on it.

\section{Numerical Examples}

\subsection{The leak test}
For the first test we set a straight filament of length 1 parametrized
by $s\in[0,1]$ as $x(s) = s, y(s) = z(s) = 0$.
We discretize with $N_n=48$ equally-spaced nodes and impose a velocity
${\bf U}=(0,1,0)$, which is orthogonal to the filament. The length of each segment is $h=1/(N_n-1)$. 
Then we solve a linear system for the force
density at the $N_n$ nodes, denoted by $\{ {\bf q}_k\}$. 
The exact solution is not a linear force
field, so there will be error in the computed velocity at points between nodes 
along the filament. The boundary velocity error
 is computed for several values of $\epsilon$ as $E_k= \|\bu({\bf X}_k)-{\bf U} \|_2$
 where $\{ {\bf X}_k \}$ is a dense set of $N_e=1505$ points on the filament.
The leak was defined as $leak =  \sqrt{ \sum_{k=1}^{N_e} E_k^2/N_e }$. In all numerical experiments
we set the viscosity to $\mu=1$.

\begin{figure}[hbtp]
\includegraphics[height=2.1in]{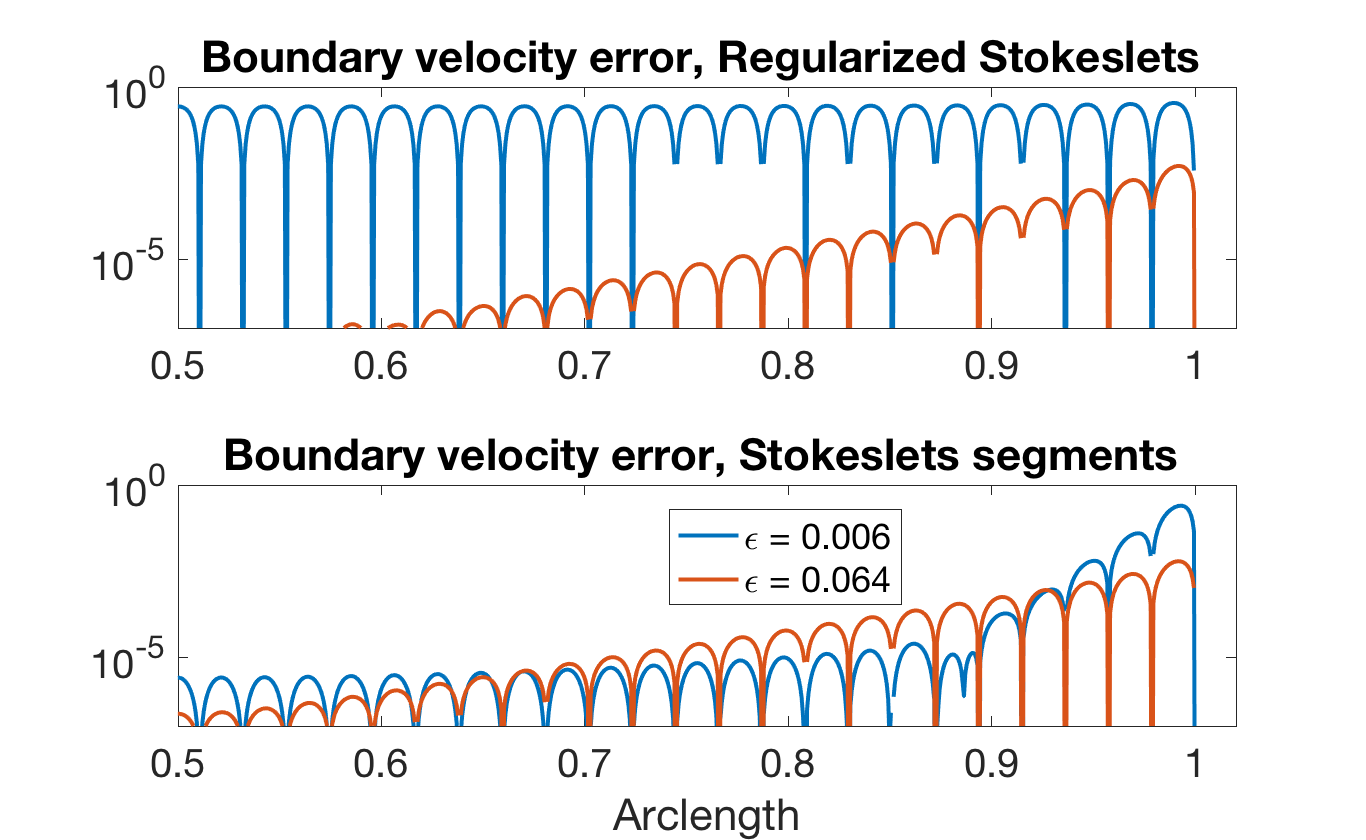}
\includegraphics[height=2.1in]{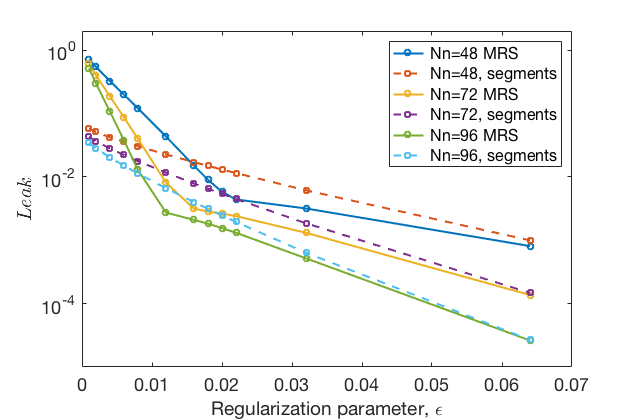}
\caption{(a) Velocity error along filament; (b) leak as a function of $\epsilon$.} 
\label{fig:leak1}
\end{figure}

{\bf The error in the filament velocity.}
When $\epsilon$ is much smaller than the internode spacing $h$, the functions $\fd$
in the Method of Regularized Stokeslet (MRS) 
do not overlap enough, producing fluid motion across the filament relative to the prescribed velocity. 
For large $\epsilon/h$,
the leak produced by the MRS decreases but the forces are spread over a larger region
surrounding the filament, which may be undesirable. In this example, the error tends to 
be larger near the endpoints of the filament, where the forces end abruptly, which is an issue
that can be addressed with nonuniform node placement. The top panel
of Figure~\ref{fig:leak1}(a) shows the error in the right half of the filament velocity using MRS for internode
spacing $h=1/47$ and two values of the regularization, $\epsilon=0.282 h$ and 
$\epsilon=3h$. The first value is so small that the error is $O(1)$ throughout the
filament. The second value is large enough that the largest error is 0.0051 near the 
filament end point and decreases quickly along the filament to a minimum of about $5.6\times10^{-8}$ at
the midpoint of the filament.

The new method using Stokeslet segments proposed in this article accounts for a continuum of regularized 
forces so there is no distance between contiguous  cutoff functions anymore. The error in the boundary
velocity (enforced for a fixed $\epsilon$) arises from the linear approximation of the forces along each segment.
The bottom panel of Figure~\ref{fig:leak1}(a) shows the error in the velocity along the filament. Even for
very small values of $\epsilon$, the error is largely confined to the regions near the end points of the
filament. 

\begin{figure}[hbtp]
\centering
\includegraphics[height=2.1in]{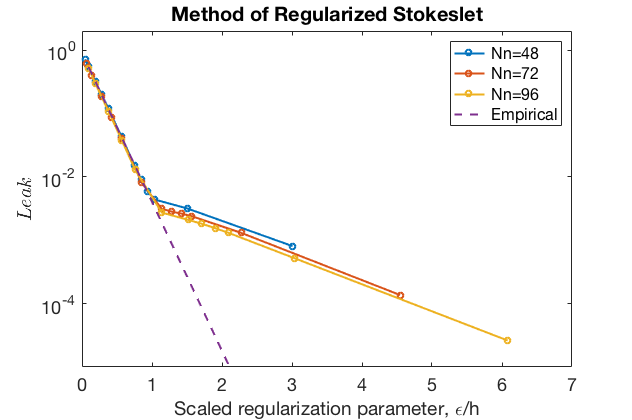}
\includegraphics[height=2.1in]{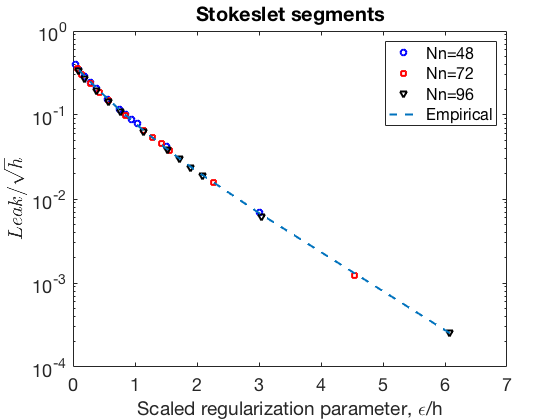}
\caption{(a) Leak using MRS; (b) Leak using Stokeslets segments.} 
\label{fig:leakEmpirical}
\end{figure}

{\bf The leak.} The leak along the filament depends on the number of nodes $N_n$ used to
enforce the velocity boundary condition. In the case of equally-spaced nodes, Figure~\ref{fig:leak1}(b) 
shows the leak for filaments discretized with 48, 72 and 96 nodes,
and for a range of values of $\epsilon$. The solid lines are the result of the MRS and are characterized
by having a large leak for $\epsilon$ near zero and decreasing exponentially for small $\epsilon$ until
a slower exponential decay leads the error for larger values of $\epsilon$. Figure~\ref{fig:leakEmpirical}(a)
shows that the fast exponential decay for different discretizations collapse onto a single line when the
regularization parameter is scaled by the internode distance $h$. This decay rate of the leak is approximated
empirically with the function
\[
 L_{MRS}(\epsilon/h)=0.9 \times 10^{-2.3\epsilon/h},\ \ \ \ \epsilon/h \in [0,1].
\]
Note that the slower decay takes place for $\epsilon/h>1$ and the leak curves collapse only approximately in
this region under the same scaling.

On the other hand, the leak using the new Stokeslet segment method collapses onto a single curve for all $\epsilon$
and all $h$ if $\epsilon$ is scaled by $h$ and the leak is scaled by the square root of $h$.
The leak in this example using the Stokeslet segments has been empirically approximated by
\[
h^{-1/2}  L_h(\epsilon)=0.25  \left( 10^{-\epsilon/h}+ 0.63\times10^{-0.46\epsilon/h} \right).
\]
The empirical curve and the collapsed data  for $N_n=48, 72$ and 96
are shown in Figure~\ref{fig:leakEmpirical}(b).

\begin{figure}[hbtp]
\centering
\includegraphics[height=2.6in]{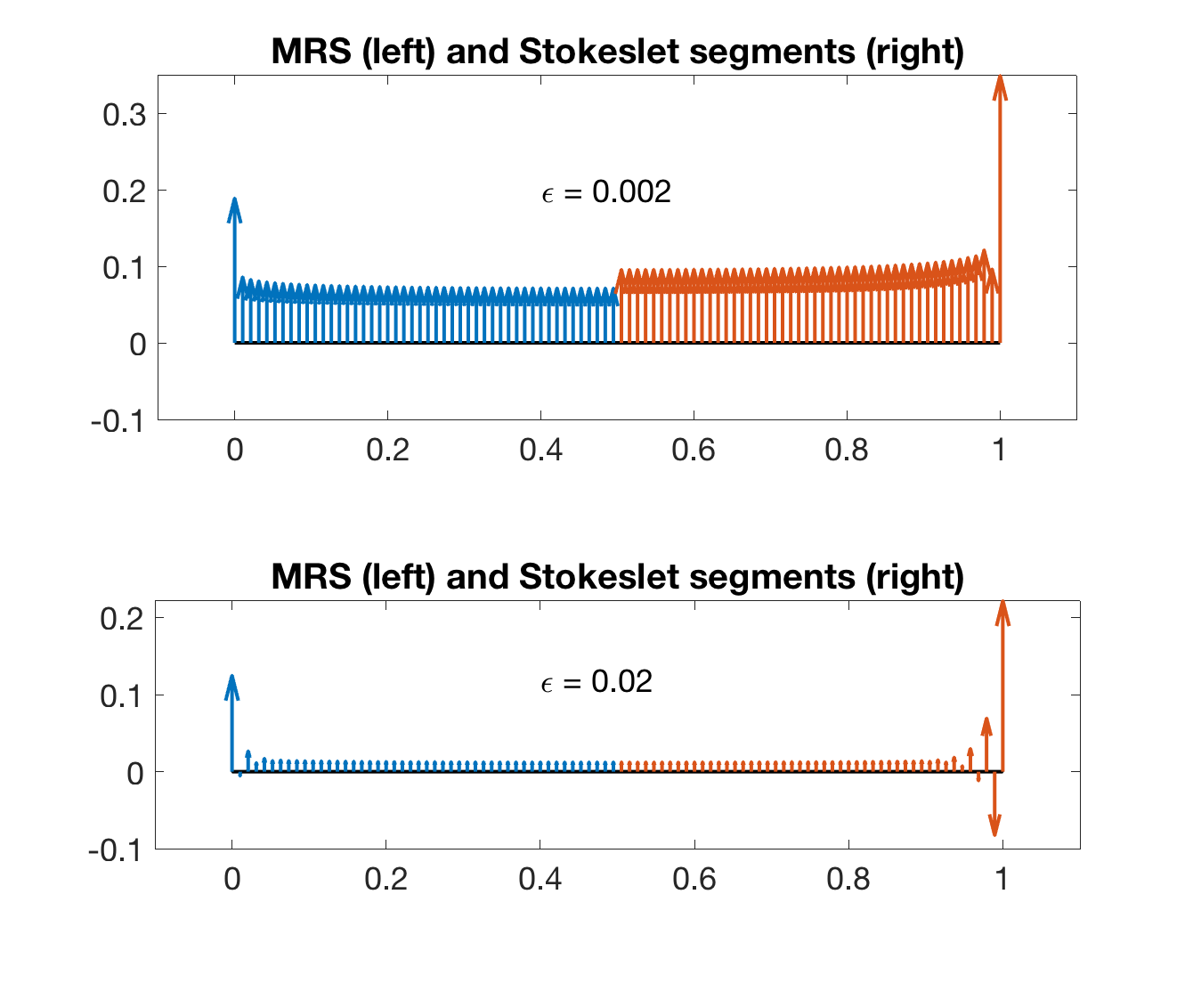}
\caption{Force density along the filament. Due to symmetry, the figure shows the results using MRS on the left
and using Stokeslets segments on the right. The  regularization parameter on the top panel is $\epsilon=0.002=0.19h$
and on the bottom panel it is $\epsilon=0.02=1.9h$. The vectors on the two panels have been scaled independently
for visualization. } 
\label{fig:forcevectors}
\end{figure}

{\bf The force density and the drag on the filament.} The force densities $\ff_k$ for $k=1,\dots,N_n$
are computed by solving a $3N_n\times3N_n$ linear system  whose matrix depends on the
position of the nodes and $\epsilon$. Figure~\ref{fig:forcevectors} shows the force densities 
obtained with the two methods and for two values of $\epsilon$. The endpoint effects are more
noticeable when using Stokeslet segments, and the force vectors alternate sign near the end points
when $\epsilon$ is ``large" enough. The MRS does the same but when the regularization is larger.
To prevent this rapid flipping of direction, $\epsilon$ must remain small enough. On the other hand,
small values of $\epsilon$ lead to larger leaks with MRS. 

\begin{figure}[hbtp]
\centering
\includegraphics[width=2.in]{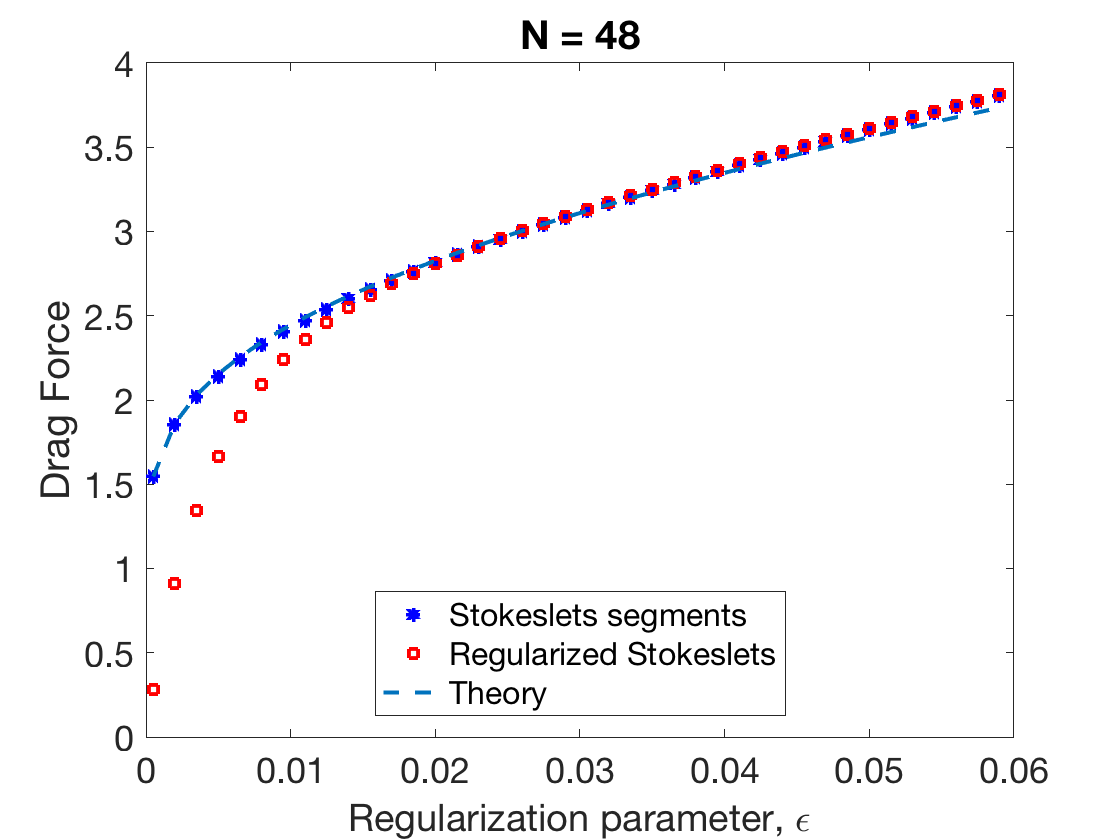}
\includegraphics[width=2.in]{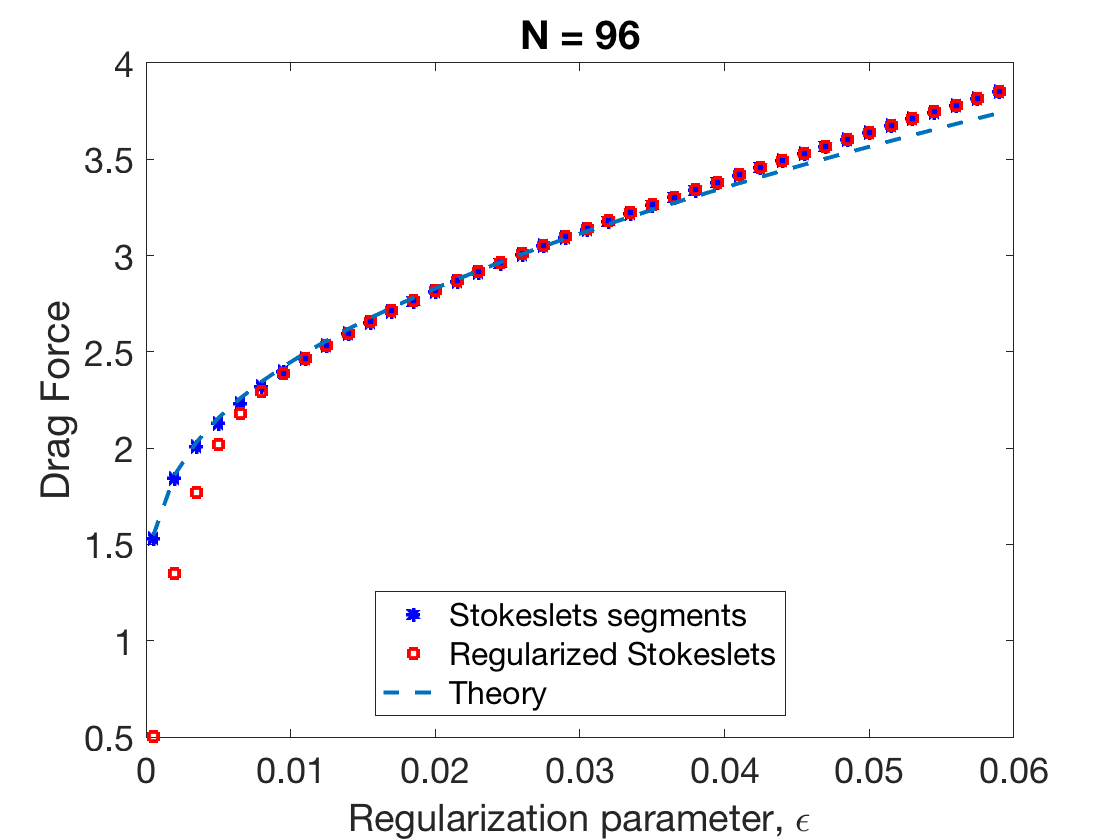}
\includegraphics[width=2.in]{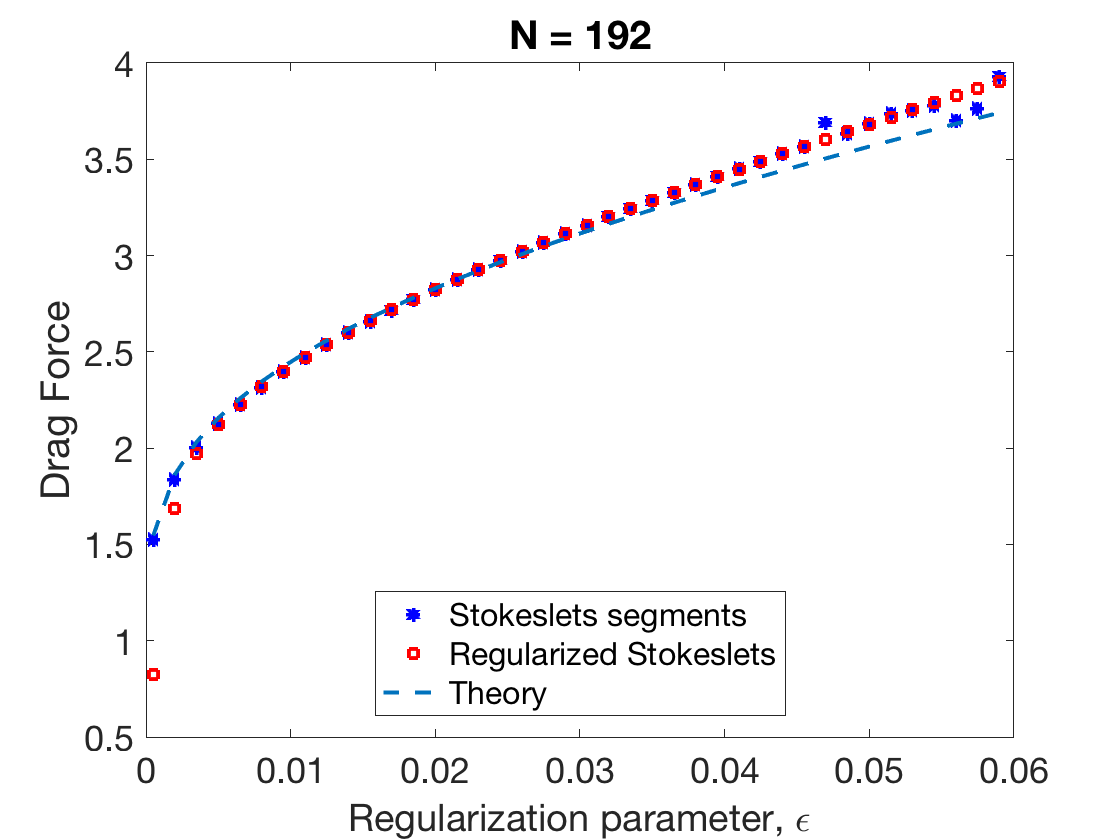}
\caption{Drag force on the straight filament using MRS and Stokeslets segments,
and comparison with the theoretical expression for the drag on a slender filament of radius $r_e=0.9634\epsilon$ ($N_n=48$), $r_e=0.9670\epsilon$ ($N_n=96$), $r_e=0.9873\epsilon$ ($N_n=192$).} 
\label{fig:dragcomparison}
\end{figure}

Conceptually, the straight filament with regularized
forces in this example may be considered a thin cylinder of length 1 and an effective
radius $r_e$ that is related to the regularization parameter $\epsilon$. The drag force
on a finite cylinder moving perpendicularly to its axis (using $\mu=1$, $\ell=1$, and $U=1$) is
\begin{equation}\label{eq:drag}
F_D = \frac{8\pi \mu \ell U}{2\log(\ell/r_e)+1} = \frac{8\pi}{1-2\log(r_e)}
\end{equation}
If we assume $r_e$ is proportional to $\epsilon$ and compute the drag force with Stokeslet segments, we find a least-squares value
of the cylinder's effective radius, we get $r_e =  0.97\epsilon$ (or $\epsilon=1.031 r_e$).
The theoretical drag using
this effective radius matches the Stokeslet segments result quite well for $\epsilon<0.04$, as
seen in Figure~\ref{fig:dragcomparison}.
The data from MRS cannot be fit well with a function of the form in Eq.~\eqref{eq:drag} for
small values of $\epsilon$ because the drag is too small. However, as the number of discretization nodes increases, the drag force from MRS approaches  the values from the Stokeslet segment method.  Similar results (not shown) apply to the filament moving parallel to its axis.

\subsection{Flagellar motion}

In many studies, a flagellum is represented by a curve where time-dependent forces and torques develop.
In simulaions that use the method of regularized Stokeslets, the parameter $\epsilon$ is sometimes regarded as the radius of a virtual slender body. Based on the
previous example, a flagellum of radius $r_e$ requires the regularization value of $\epsilon = 1.031 r_e$.
Given that
in the MRS the distance $\Delta s$ between nodes used to discretize the filament must be on the order of
$\epsilon$, the choice of regularization imposes a minimum number of nodes needed for the simulation.
A typical MRS simulation of a flagellum of length $\ell$ uses 100-200 nodes with node separation $\Delta s=0.005 \ell-0.01 \ell$ and $\epsilon\approx 1.3-5 \Delta s$~\cite{SIMONS20151639, Ko2017}. This gives a length-to-radius ratio of
about 50-100, which is too small for some applications.

In this example, we model the flagellum as an inextensible elastic filament based on the formulation described
in~\cite{SIMONS20151639} so we only mention the main equations here. If $\bX(s,t)$ is a point on the filament with $s$ denoting the arc length parameter, we prescribe a target curvature $\Omega(s,t)$. 
A planar filament $\bX(s,t)=(x(s,t),y(s,t),0)$ of length $\ell$ discretized with $M-1$ segments of length $h$
($M$ nodes) and we define the energy
\[
E = \frac{1}{2} \sum_{j=1}^{M-1}  \kappa_T \left(\| D^+\bx_j \|-1 \right)^2 h
+  \frac{1}{2} \sum_{j=1}^{M-1}  \kappa_B \left( D^2 y_j D^0 x_j - D^2 x_j D^0 y_j -\Omega(s_j,t) \right)^2 h
\]
where $D^+\bx_j = (\bx_{j+1}-\bx_j)/h$, $D^0\bx_j = (\bx_{j+1}-\bx_{j-1})/2h$,  
$D^2\bx_j = (\bx_{j+1}-2\bx_j+\bx_{j-1})/h^2$, and $\kappa_T$ and $\kappa_B$ are tensile and bending stiffnesses.
Complete details are found in~\cite{SIMONS20151639}.

The prescribed dimensionless target curvature is 
\begin{equation}\label{eq:omega}
\Omega(s,t) = \frac{A k^2\sin(k s-2\pi t)}{\sqrt{1- A^2 k^2 \cos^2(ks-2\pi t)}} + \Omega_0
\end{equation}
for $0\le s\le 1$, which is consistent with a sinusoidal traveling wave. The discrepancy between the curvature of the filament and the target curvature
produces a penalty force $\FF_k = -\partial_{\bx_k}  E$. We define the corresponding force density as
$\ff_k = \FF_k/h$ except at the end points where it is $\ff_1 = 2\FF_1/h$ and $\ff_M = 2\FF_M/h$.

\begin{table}[hbtp]
\centering
\begin{tabular}{|r|c|l|}
\hline
   Parameter &  Value &  Description \\
\hline
    $\kappa_B$  & $0.0221$   &   Bending stiffness\\
    $\kappa_T$  & $2.950$  &   Tensile stiffness \\
    $A$  & $0.075$  &    Curvature amplitude \\
    $k$  & $9\pi/4$  &   Wave number \\
    $\sigma$  & $2\pi$  &    Beat frequency  \\
    $\ell$  & $1$  &  Flagellum length   \\
    $\Delta t$  & $2.5\times 10^{-7}$  &  Time step  \\
    $t_{f}$  & $70$  &  Final time   \\
    $M$  & $24$  &  Number of discretization nodes   \\
    $\epsilon$  & $1/300$  &  Regularization parameter   \\
    $\mu$  & 1 & Fluid viscosity  \\
\hline
\hline
\end{tabular}
\caption{Dimensionless parameter values used in the planar beat numerical examples.}
\label{table:params1}
\end{table}

\begin{figure}[h]
\hskip20pt\begin{minipage}[t]{0.5\textwidth}
{ \includegraphics[width=3.2in]{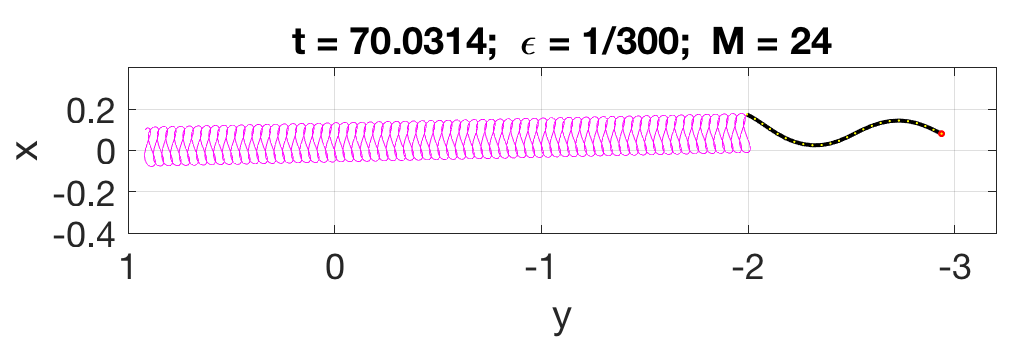} }\\
{ \includegraphics[width=3.2in]{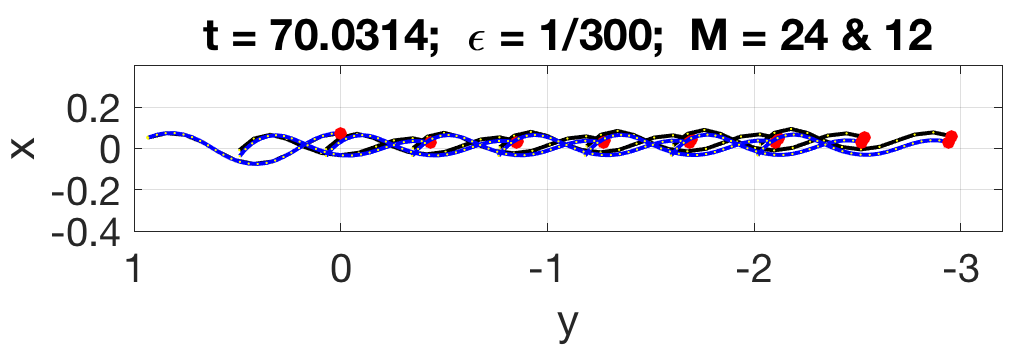} }
\end{minipage}
\begin{minipage}[t]{0.5\textwidth}
{\raisebox{-0.45\height}{{ \includegraphics[width=2.3in]{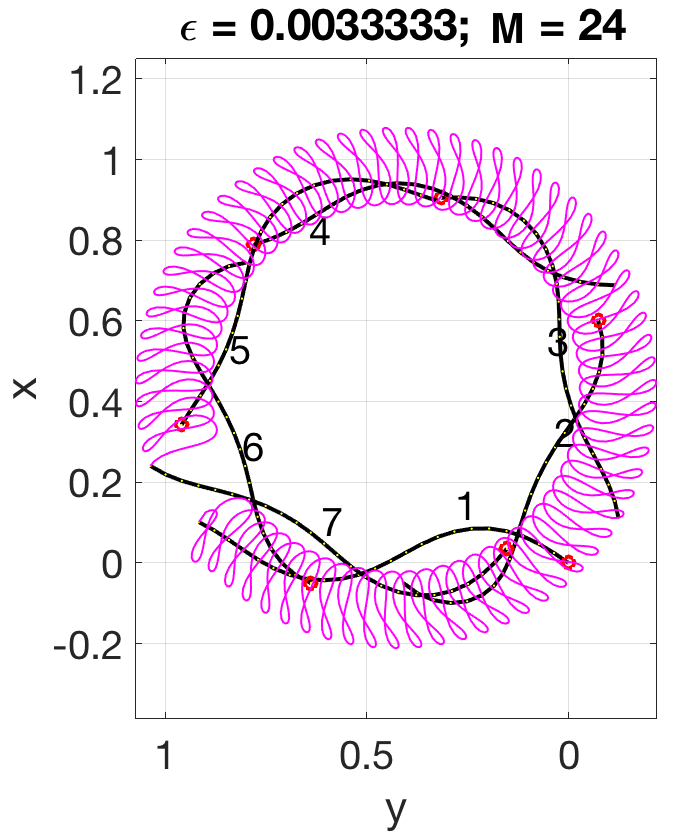} }}}
\end{minipage}
\caption{Straight line trajectory of a flagellum driven by the prescribed curvature in Eq.~\eqref{eq:omega}
with $\Omega_0=0$ using $M=24$ nodes (top left). A trajectory comparison between filaments discretized
using $M=24$ and $M=12$ nodes (bottom left).  Circular trajectory when $\Omega_0=1.4989$ (right).} 
\label{fig:runandturn}
\end{figure}

Using $M=24$ we initialized the flagellum in the $xy$-plane to have the curvature in 
Eq.~\eqref{eq:omega} with $\Omega_0=0$.
The regularization parameter was set to $\epsilon = \ell/300$. This corresponds to
a slender cylinder of length  to radius ratio of $\ell/r_e = 1.031(300) = 309$.
Some sperm are known to swim both in straight lines and in circular motion~\cite{Iida2017ACS}
depending on the amount of asymmetry in the curvature of the flagellum~\cite{RIKMENSPOEL1985395}.
We set out to use our method to produce both behaviors. The straight runs are accomplished by prescribing only the sinusoidal term in the target curvature Eq.~\eqref{eq:omega}, which produces traveling waves along the flagellum resulting in  straight swimming in the direction opposite the traveling wave. Figure~\ref{fig:runandturn} (top left) shows a snapshot at the end of the simulation and the oscillatory trajectory of the last point of the flagellum. In this example it takes about 22 flagellar beats to swim one body length. The bottom left panel shows a comparison of the
flagellum trajectories when the filament is discretized using $M=24$ and $M=12$ nodes while keeping all other parameters fixed. In both simulations the swimming speed is the same; the only difference is a very slight angle in the trajectory, which is influenced by the discretization and the initial condition.

\begin{figure}[hbtp]
\centering
\includegraphics[width=3in]{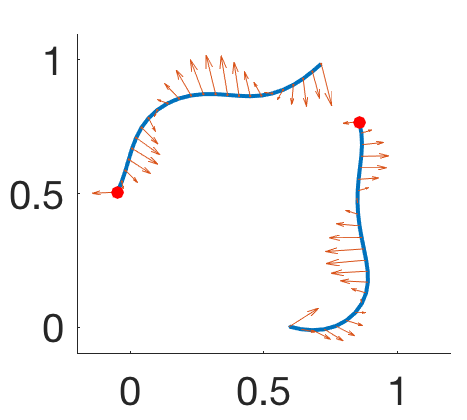}
\caption{Forces at $t=33.14$ and $t=58.68$ of the circular trajectory.} 
\label{fig:circleforces}
\end{figure}

Our model for introducing an asymmetry that will induce circular motion is to add
 a constant term to the curvature of every point. This is appropriate for a sperm flagellum since it contains active dyneins capable of producing bending throughout its length.
 Figure~\ref{fig:runandturn} (right) shows seven snapshots of the flagellum moving in a circle. The snapshots are shown every 11.25 flagellar beats and the added curvature is $\Omega_0=0.4Ak^2=1.4989$.
Figure~\ref{fig:circleforces} shows the forces that develop along the flagellum at two different times during in the
simulation.

\subsection{The Stokeslet segment method for flows in the half-space $z>0$}
As described in~\cite{Ainley:08, ImagesGeneral}, the flow bounded by a plane can be realized using a regularized image system involves a Stokeslet, a Stokes doublet, a potential dipole and rotlets.  We follow the same procedure
and assume the segment $\by^*(\alpha)= \by^*_0 + \alpha(\by^*_1-\by^*_0)$ is in the fluid domain $z>0$ and set the notation $\by^*_0\cdot\hat{\bf e}_3=H_0$ and $\by^*_1\cdot\hat{\bf e}_3=H_1$. The
image of the segment is $\by(\alpha) = \by_0 + \alpha(\by_1-\by_0)$, where
$\by_0=\by^*_0-2H_0\hat{\bf e}_3$ and $\by_1=\by^*_1-2H_1\hat{\bf e}_3$.

\begin{picture}(200,100)(-10,-10)
\put(20,15){\circle*{5}}
\put(60,5){\circle*{5}}
\put(90,-2){\circle*{5}}
   \put(0,30){\line(1,0){110}}\put(120,25){$z=0$}
   \put(120,50){fluid domain}
   \put(120, 0){images}
\put(20,45){\circle*{5}}
\put(60,55){\circle*{5}}
\put(90,62){\circle*{5}}
  \put(60,10){\vector(-1,4){8}}  
\put(13,5){$\by_0$}    \put(50,-5){$\by(\alpha)$}      \put(83,7){$\by_1$}
\put(13,55){$\by^*_0$}    \put(50,65){$\by^*(\alpha)$}      \put(87,52){$\by^*_1$}
  \put(51,45){\circle*{3}}   \put(40,40){$\hat{\bx}$} 
   \put(20,15){\line(4,-1){70}}
   \put(20,45){\line(4,1){70}}
  \end{picture}

\begin{figure}[hbtp]
\centering
\includegraphics[width=3in]{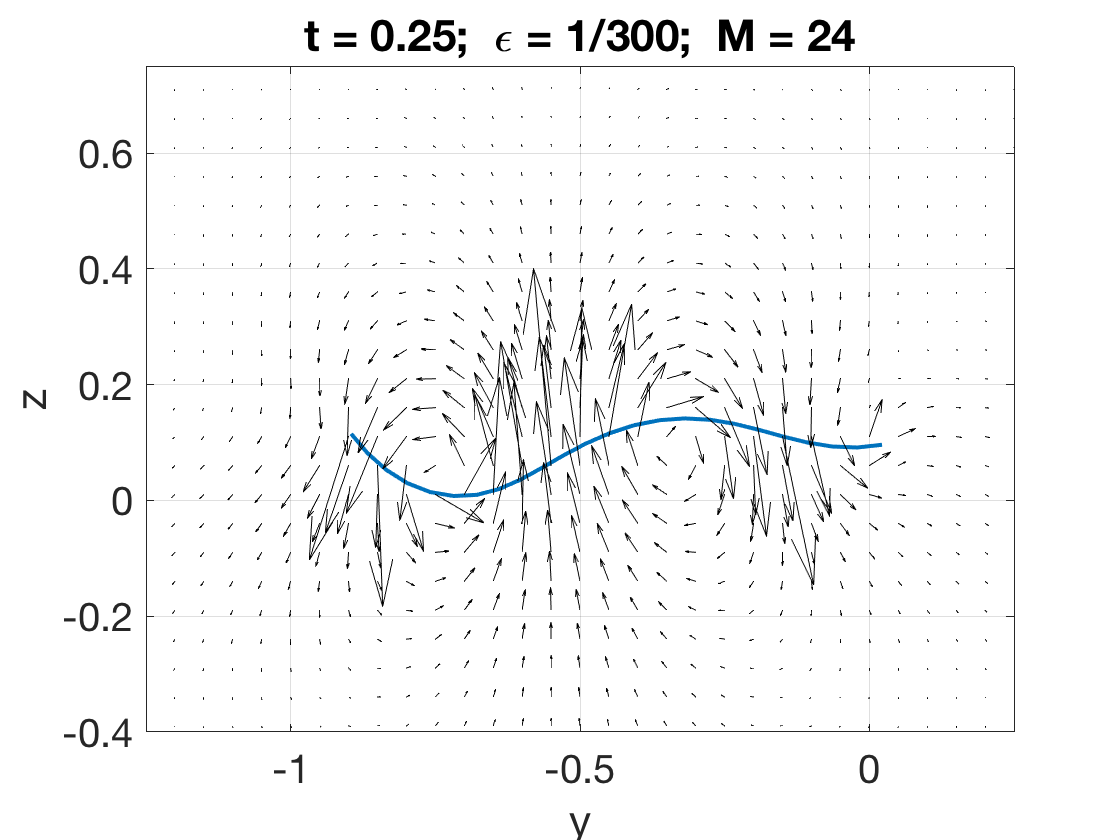}
\includegraphics[width=3in]{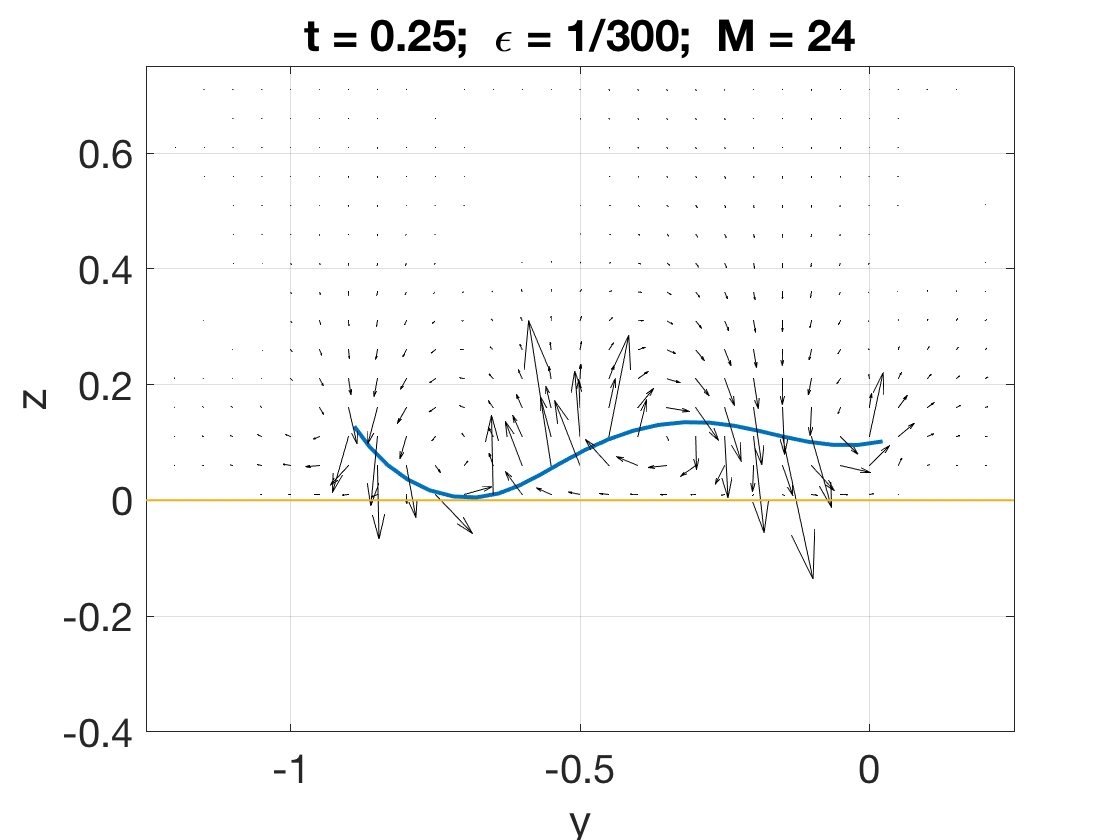}
\caption{Comparison of the velocity field around the flagellum in the absence of boundaries and flow bounded by a plane wall.} 
\label{fig:gridvelcomp}
\end{figure}

Using the notation $8\pi u_i = S^\epsilon_{ij}f_j$ to represent the $i$-th component of the Stokeslet, the Stokes doublet is defined as
\[
\Delta_{ijk} = \frac{\partial S^\epsilon_{ij}}{\partial x_k} 
= \frac{1}{R^3}\left( x_i \delta_{kj} -x_k\delta_{ij} + x_j \delta_{ki} \right)
- \frac{3}{R^5} x_k  \left( \epsilon^2 \delta_{ij} +  x_i x_j   \right)
\]
and is applied to the forcing vector ${\bf q} = (q_1,a_2,q_3)=(-f_1, -f_2, f_3)$.
The complete system of images is
\[
-S^\epsilon_{ij} f_j + 2 H \Delta_{i3j} q_j + H^2 PD_{ij} q_j + 2 H R_{ij} f_j
\]
where the dipole $PD_{ij}$ is the integrand of Eq.~\eqref{def:dipole} and $R_{ij}$ represents
 the difference between two regularized rotlets~\cite{Ainley:08, ImagesGeneral}
\[
R_{ij} f_j =  \frac{3\epsilon^2}{R^5} (   x_3 f_1 \delta_{1i} + x_3 f_2 \delta_{2i}  -(f_1x_1+f_2x_2) \delta_{3i}  ) .
\]
In the images formula, $H=\by(\alpha)\cdot\hat{\bf e}_3$.
The fluid velocity of this combination of elements exactly cancels at the wall the flow of the original Stokeslet.
We note that the elements involved are of the form $P_1(\alpha) R(\alpha)^{-1}$, $P_2(\alpha) R(\alpha)^{-3}$,
and $P_3(\alpha) R(\alpha)^{-5}$, where $P_1$, $P_2$ and $P_3$ are polynomials of degree up to 5 since
$x_i$, $f_i$ and $H$ are linear functions of $\alpha$. As before, terms of these forms are computed using the
recursion in Eq.~\eqref{Recurrence}. Figure~\ref{fig:gridvelcomp} shows the velocity field on the plane of motion of a flagellum in open space compared to one swimming near a wall. The velocity field has the same scaling in both panels.

\begin{figure}[hbtp]
\centering
\includegraphics[width=2in]{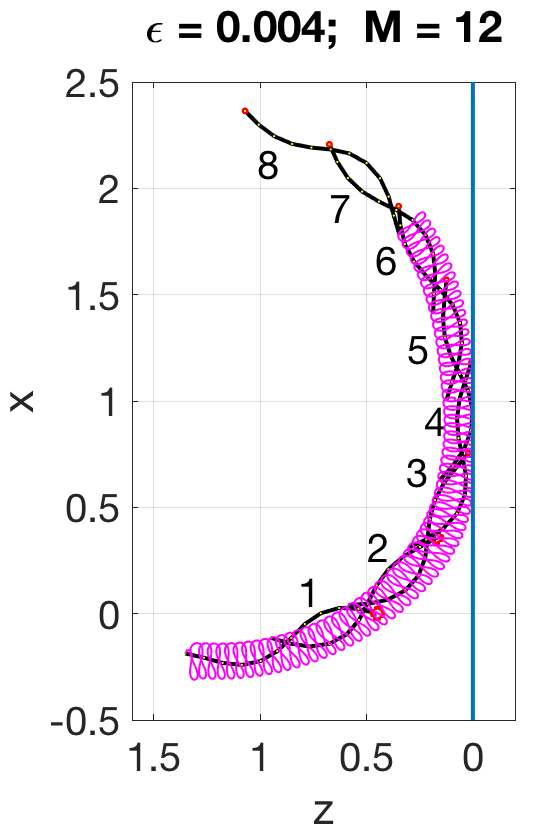}
\caption{Trajectory of a sperm driven by the prescribed curvature in Eq.~\eqref{eq:omega}
with $\Omega_0=0.6$ and bounded by a plane. The parameters used are hose in 
Table~\ref{table:params1} except $M=12$ and $\epsilon=0.004$.} 
\label{fig:runandturnwall}
\end{figure}

Figure~\ref{fig:runandturnwall} shows the trajectory of a flagellum moving in a circular
path and approaching a solid wall. In this example we use only 11 segments (discretizing a
curve of unit length) and a value of $\epsilon=1/250$ to demonstrate that the proposed method
can be implemented even in an extreme case of
using very few segments. The waveform of the flagellum is affected by the presence
of the wall as the flagellum moves parallel to the wall. Eventually, the sperm swims away from
the wall.

\clearpage

\subsubsection{Nonplanar flagellar beats and turning: the Kirchhoff rod formulation}
It is known that the trajectories of swimming sperm involve straight stretches characterized
by symmetric (sinusoidal) wave forms, looping turns during which
the wave form of the flagellum is asymmetric, and sharp turns~\cite{Kaupp2016, Wood725, RIKMENSPOEL1985395}.
In order to model a flagellum with target curvatures in all coordinate directions, we use
the Kirchhoff rod formulation described
in~\cite{OLSON2013169} so we only mention the main equations here. If $\bX(s,t)$ is a point on the filament with $s$ denoting the arc length parameter, a local
orthonormal coordinate system $\{ \bD_1,\bD_2,\bD_3\}$ is defined centered at $\bX$. At $t=0$, the
unit vectors are defined as $\bD_3 = \partial\bX(s,0)/\partial s$, $\bD_1$ and $\bD_2$ are perpendicular
to the axis of the rod. The target curvature is a vector written in the local coordinates 
$\Omega = \Omega_1 \bD_1 +\Omega_2 \bD_2 +\Omega_3 \bD_3 $.

\begin{figure}[hbtp]
\centering
 \includegraphics[width=2in,angle=0]{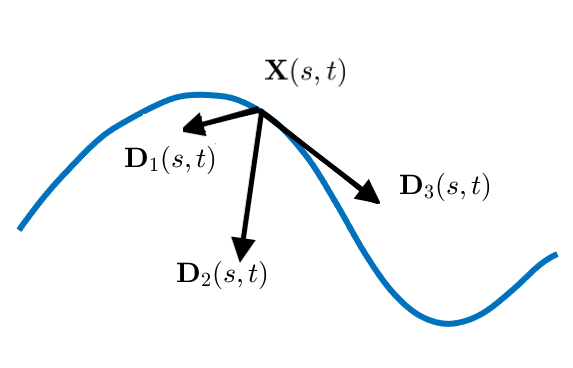}
\caption{Schematic of the local coordinate system in the Kirchhoff rod model.} 
\label{fig:KRschematic}
\end{figure}

If $\ff(s,t)$ and $\tau(s,t)$ are the force density and torque density exerted by the
fluid on the rod, then the conservation of force and torque become
\begin{equation}\label{fConserv}
0 = \ff + \frac{\partial \FF}{\partial s},\ \ \ \ 
0 = \tau + \frac{\partial \TT}{\partial s} + \left( \frac{\partial\bX}{\partial s}\times\FF \right)
\end{equation}
where $\FF$ and $\TT$ are internal forces and couples. Writing variables in terms of the local coordinates
\[
\FF = \sum_{i=1}^3 F_i \bD_i,\ \ \ \TT = \sum_{i=1}^3 Q_i\bD_i,\ \ \ 
\ff = \sum_{i=1}^3 f_i \bD_i,\ \ \ \tau = \sum_{i=1}^3 q_i\bD_i,\ \ \ 
\]
the constitutive relations are
\begin{eqnarray}
Q_1 &=& a_1 \left( \frac{\partial \bD_2}{\partial s}\cdot\bD_3 - \Omega_1 \right), \ 
Q_2 = a_2 \left( \frac{\partial \bD_3}{\partial s}\cdot\bD_1 - \Omega_2 \right),\ 
Q_3 = a_3 \left( \frac{\partial \bD_1}{\partial s}\cdot\bD_2 - \Omega_3 \right),\label{Qi}\\
F_i &=& b_i \left( \frac{\partial \bX}{\partial s}\cdot\bD_i - \delta_{3i} \right),\ \mbox{ for } i=1,2,3.\label{Fi}
\end{eqnarray}
where the $a_i's$ and $b_i's$ are stiffness constants and $\Omega_i$ are prescribed curvature 
functions that initiate and maintain the motion of the rod.

\begin{table}
\centering
\begin{tabular}{|r|c|c|c|}
\hline
   Parameter &  Value & Units & Description \\
\hline
    $a_1 = a_2 = a_3$  & $4.9587$   & mg $(\mu$m$)^3/s^2$   &  Bending stiffness\\
    $b_1 = b_2 = b_3$  & $0.8264$  &  mg $(\mu$m$)/s^2$   & Shear stiffness \\
    $A$  & $3.5$  &  $1/\mu$m  &  Curvature amplitude \\
    $k$  & $9\pi/160$  &  $1/\mu$m  &  Wave number \\
    $\sigma$  & $550$  &  $1/s$  &  Beat frequency  \\
    $\ell$  & $40$  & $\mu$m  & Length   \\
\hline
\hline
\end{tabular}
\caption{Dimensional parameter values used in the Kirchhoff rod numerical example. 
For the computations the parameters our non-dimensionalized using $\ell$ as the characteristic length scale,
$T_0=2\pi/\sigma$ as the characteristic time scale, and $\mu \ell^2/T_0$ as the reference force.}
\label{table:params}
\end{table}

 A filament of length $\ell$ is discretized with $N-1$ segments ($N$ nodes including the
 end points).  For $k=1,2,\dots,N-1$ and given the current configuration of the curve, we first compute
$ F_i(s_{k+1/2})$ and $Q_i(s_{k+1/2})$ with finite difference approximations of
Eq.~\eqref{Qi}-\eqref{Fi}. Then $\ff(s_k)$ and $\tau(s_k)$ are computed with finite difference 
approximations of Eq.~\eqref{fConserv}. Finally the velocity field at any point due to the
force and torque densities along the rod is evaluated by adding the flow contributions of
Stokeslets and rotlets in each segment. 
To compute the rate of change of the local coordinate vectors $\{ \bD_1,\bD_2,\bD_3\}$,
the curl of the fluid velocity formula is needed, leading to a distribution of rotlets and
dipoles along the rod.
More details about the implementation can be found in~\cite{OLSON2013169}.

In the Kirchhoff rod framework, this asymmetry
can be modeled by augmenting he planar case with a target curvature in the orthogonal direction. This adds
a bend to every point on the flagellum.
The twist in the filament
is $\Omega_3(s,t)$ and is kept at zero and the other components are set to
 \begin{eqnarray*}
\Omega_1(s,t) &=& W_1  \\
\Omega_2(s,t) &=& \frac{A k^2\sin(k s-\sigma t)}{\sqrt{1- A^2 k^2 \cos^2(ks-\sigma t)}} + W_2
\end{eqnarray*}
where $W_1$ and $W_2$ are random piecewise constant functions of time with values
chosen every 15 beats from a uniform distribution in $[-\Omega_0,\Omega_0]$ where 
$\Omega_0=0.4 A k^2=0.04372\ 1/\mu$m.

 \begin{figure}[hbtp]
\centering
\includegraphics[width=4in]{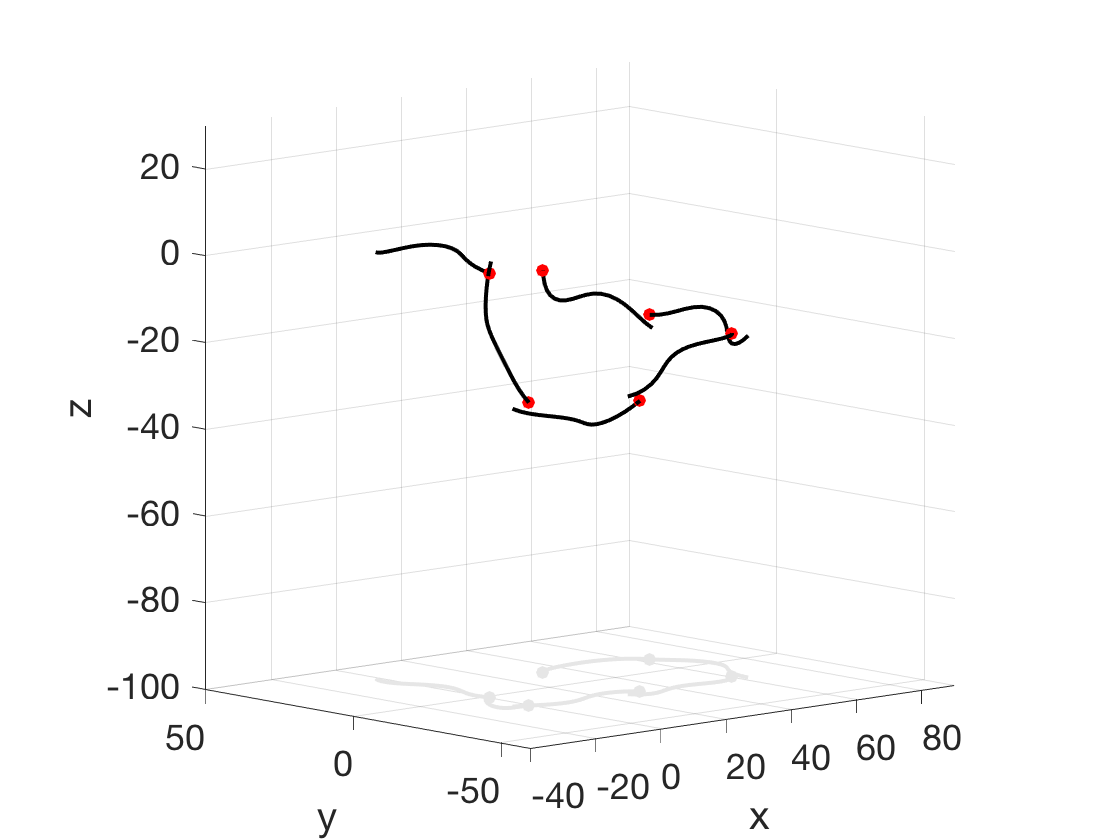}
\caption{Snapshots of the trajectory of the flagellum as it samples three-dimensional
space using the Kirchhoff rod model. The shadows at the bottom are projections of
the flagellum shape on the plane $z=-100$ for visualization.} 
\label{fig:KRsnapshots}
\end{figure}

We discretized the
filament with $N=20$ nodes and set the regularization parameter to $\epsilon=0.005$
which corresponds to a flagellar effective radius of
$r_e = 0.97\epsilon = 0.00485$ for a length-to-radius ratio of $1/r_e=206$. Figure~\ref{fig:KRsnapshots}
shows six snapshots approximately every 16.25 beats. The initial shape lies in 
the plane spanned by $\bD_2$ and $\bD_3$ at $t=0$. The curvature component $\Omega_2$ 
drives the (local) planar motion while $\Omega_1$ introduces out-of-plane motion.
%

\subsection{Discussion and conclusions}
In the method of regularized Stokeslets, the regularization parameter $\epsilon$ and the node separation $\Delta s$ must be chosen in such a way that the cutoff functions overlap sufficiently. We can think of these cutoffs as basis functions centered at the discretization nodes so that there must be overlap in order to resolve the force field in between nodes. The leak test revealed that when $\epsilon>\Delta s$, the error in the velocity between nodes is of the form $Leak \sim 10^{-c_1 \epsilon/\Delta s}$, which decreases as the numerical parameters decrease if we choose $\epsilon = c_2 \Delta s^p$, with $0<p<1$. This way, as $\Delta s\to0$, $\epsilon/\Delta s\to\infty$ and the leak would vanish. 

The method proposed here represents the limit as the separation between cutoff functions goes to zero on straight segments. There is no discrete separation between the functions since $\Delta s$ has been removed. The new discretization parameter $h$, which breaks up a curve into a series of line segments of length $h$ with linear forcing, can be selected based on the curvature of the filament and the variation of the force field along it, but otherwise independently from the regularization size $\epsilon$. The parameter $h$ is the length of segments on which the forces can be accurately represented as piecewise linear. This is a major difference in the two methods since the 
MRS requires a discretization size  $\Delta s\approx \epsilon$, leading to a number of discretization nodes an order of magnitude larger than the new proposed method.

Figure~\ref{fig:runandturn} (bottom left) supports this conclusion. Fixing $\epsilon$ sets the radius of a slender
cylinder and varying the number of segments from 24 to 12 has an insignificant effect on the flagellum swimming speed and waveform amplitude. This is because both discretizations reasonably resolve the forces that develop on the flagellum.  The only difference is in the angle of the axis of motion, which we measured  to be  $\Delta\theta\approx\pi/320$. The angle is selected during a brief initial transient period before the flagellum settles
to a steady motion and has a numerical dependence on the initial shape and the discretization.

Further evidence that the segment size $h$ is not strongly related to $\epsilon$ is in the simulations of turning sperm (Figure~\ref{fig:circleforces}) where a sufficient number of segments to resolve the variations in the forces per wavelength may be on the order of 20-30. More importantly, the number of segments can be chosen independently from $\epsilon$, although the simulations produce better results when $\epsilon$ is substantially smaller than the segment length $h$. This is the case in Figure~\ref{fig:forcevectors} where the resulting force vectors began to change sign  when  $\epsilon>h$ and in Figure~\ref{fig:dragcomparison} where the drag computation became erratic when $\epsilon>8h$.  Both of these examples required inverting a linear system for the forces in order to impose velocity conditions on the filament. The issue of forces developing oscillations were observed only in cases when the matrix of the linear system was badly conditioned. Further investigation is necessary to characterize the cases where the relative sizes of $h$ and $\epsilon$ lead to ill conditioned systems; however, the choice $\epsilon\ll h$ seems to work well. 

We have presented a method of Stokeslet segments under the assumption that the force field along a
filament is piecewise linear. If desired, extensions to piecewise quadratic or higher degree polynomial are straight forward.  The presentation focused on the velocity field; however, the pressure can be treated in the same way.  The pressure corresponding to the regularized Stokeslet, rotlet, and dipole in Eq.~\eqref{regStokeslet}-\eqref{regDipole1}
are
\[
p_S = (\ff\cdot\bx) \frac{2R_0^2+3\epsilon^2}{8\pi R_0^7},\ \ \ 
p_R = 0,\ \ \ 
p_D = -(\ff\cdot\bx) \frac{105\epsilon^4}{8\pi R_0^9}
\]
so that their line integrals can be computed with the recurrence relation~\eqref{Recurrence}.

\section*{Acknowledgement}
The author acknowledges partial support from NSF award DMS-1043626.

\end{document}